\documentclass[letterpaper,11pt]{article}

\usepackage[utf8]{inputenc}
\usepackage[T1]{fontenc}
\usepackage{dsfont}
\usepackage{textcomp}
\usepackage{amsmath}
\usepackage{bussproofs}
\usepackage{url}
\usepackage[dutch,french,german,polutonikogreek,english]{babel}
\usepackage[multiple]{footmisc}
\usepackage{latexsym}
\usepackage[authoryear,square]{natbib}
\usepackage{filecontents}

\newcommand{\dutch}[1]{\foreignlanguage{dutch}{#1}}
\newcommand{\french}[1]{\foreignlanguage{french}{#1}}
\newcommand{\german}[1]{\foreignlanguage{german}{#1}}
\newcommand{\weg}[1]{}

\DeclareUnicodeCharacter{03BB}{$\lambda$}
\DeclareUnicodeCharacter{03C9}{$\omega$}
\DeclareUnicodeCharacter{03C0}{$\pi$}

\newdimen\savedxvi
\newcommand{\rhop}{%
\savedxvi=\fontdimen16\textfont2
\mbox{$\fontdimen16\textfont2=5pt\rho_{P}$}%
\fontdimen16\textfont2=\savedxvi
}

\title{Intuition, iteration, induction\\ (draft)}
\author{Mark van Atten\thanks{SND (CNRS/Paris IV), 28 rue
Serpente,
75006 Paris,
France.
\texttt{vanattenmark@gmail.com}.}}

\begin{document}
\maketitle

\begin{abstract}
In
\emph{Mathematical Thought and Its Objects},
Charles Parsons argues that our knowledge of
the iterability of functions on the natural numbers
and of the validity of complete induction
is
not intuitive knowledge;
Brouwer disagrees on both counts.
I will compare Parsons'
argument with Brouwer's and defend the latter.
I will not
argue that Parsons is wrong once his own conception of
intuition is granted,
as I do not think that that is the case.
But I will try to make two points:
(1) Using elements from Husserl and from Brouwer,
Brouwer's claims can be justified in more detail than he
has done;
(2) There are certain elements in Parsons'
discussion that,
when developed further,
would lead to
Brouwer's notion thus analysed,
or at least something
relevantly similar to it.
(This version contains a postscript of May, 2015.)
\end{abstract}

\tableofcontents

\newpage

\section{Introduction}

Consider the following two principles:
\begin{enumerate}
\item Iterability:
Any total function $f:\mathds{N}\rightarrow \mathds{N}$
can be iterated arbitrarily many times;

\item Complete induction:
\begin{prooftree}
\alwaysNoLine 
\AxiomC{$A(0)$} 
\AxiomC{[$A(a)$]}
\UnaryInfC{$\vphantom{\int^0}\smash[t]{\vdots}$}
\UnaryInfC{$A(Sa)$} 
\AxiomC{$N(t)$} 
\alwaysSingleLine
\TrinaryInfC{$A(t)$}
\end{prooftree}
where
$A$ is a placeholder for any well-defined predicate
on the natural numbers,
and $N$ is a natural number-predicate.%
\footnote{%
Instead of a rule,
induction may also be formulated as an axiom schema
(using implication).
In the Brouwerian context,
however,
the latter also has to be construed as a rule:
As mathematical objects and proofs exist only as
the result of the subject's activity,
to hypothesise that a proof exists therefore
must be
to hypothesise that
the subject
\emph{has constructed}
a proof.
Making an assumption
towards proving an implication
then has epistemic import.
As pointed out in \citet{Sundholm.Atten2008},
this means that in this respect Brouwerian logic
differs from Gentzen's
Natural Deduction,
in which this is not the case.}
\end{enumerate}
In his book
\emph{Mathematical Thought and Its Objects} 
\citep{Parsons2008},
Charles Parsons argues that
our knowledge that these two general principles
are valid
is not intuitive knowledge.
Brouwer,
on the other hand,
in his dissertation 
\citep{Brouwer1907}
claims that they are.

I will compare Parsons' position
on iteration and induction
with Brouwer's
and defend Brouwer's view,
in the following sense.
Brouwer's and Parsons'
conceptions of intuition are,
of course,
different.
I will not argue that Parsons is wrong
once his own conception of intuition is granted,
as I do not think that that is the case.
But I will try to make two points.
The first is that,
using elements from Husserl and from Brouwer,
Brouwer's claims that
the principle of induction and the iterability of functions
are intuitive
can be justified in more detail than he
has done.
The second is that
there are certain elements in Parsons' discussion
that,
when developed further,
would lead to Brouwer's notion thus analysed,
or at least something similar to it.

\section{Parsons}
\label{Parsons}
For Parsons,
intuition plays two rôles in the foundations of mathematics.
Intuition of objects falling under a given concept
shows that the mathematical concept in question is not
empty;
and intuitive knowledge is more evident than knowledge of
other types,
in particular knowledge whose justification involves
appeals to principles of reason \citep[pp.~113, 336]{Parsons2008}.
But unlike Brouwer,
Parsons is not a constructivist,
and intuition has no overall legislative rôle for him.
He does not generally require that,
for a mathematical concept to be instantiated,
it should be
instantiable
by an
intuited object;
this is only required for concepts
at the bottom of our conceptual edifice.
Other concepts then arise as combinations and idealisations
of the lower ones,
but there is no requirement that we have intuitions
of objects that fall under these higher concepts.
Naturally,
an account is then needed of
which idealizations are legitimate and which ones go too
far.
Parsons says this will depend on a theory of reason.
He develops a number of ideas on this in the last chapter of
\emph{Mathematical Thought and Its Objects},
`Reason',
but my concern here will rather be with his views on
the intuitive part of mathematics.

Parsons distinguishes intuition
\emph{of} objects
from intuition
\emph{that}
a proposition is true.
For mathematical intuition specifically,
Parsons takes the objects of intuition to be strings of
strokes,
e.g.,
$||||||$.
He allows for the possibility to see specific inscriptions
of such
strings
(tokens)
\emph{as}
types,
following Husserl here
in holding that sometimes
\emph{imagination}
of the token
can found intuition of a type \citep[p.~173]{Parsons2008}.

The first thing Parsons says in defining intuitive knowledge
is
that we have intuitive knowledge
\emph{that}
$p$ if $p$ can be `seen' 
(quotation marks Parsons')
to be true
on the basis of intuiting objects that it is about
\citep[p.~171]{Parsons2008}.
He then generalises this by counting as intuition of an
object
not only an actual perception of one,
but also an imagination of an arbitrary one.
Otherwise,
we could not justify the claim that
it is intuitive knowledge that
any such string can be extended
by placing one more stroke at the right.
An arbitrary string,
Parsons says,
may be imagined in the way we imagine
a large crowd at a baseball game;
it need not be part of the content of that imagination
that the crowd consists of exactly $n$ people \citep[pp.~173–174]{Parsons2008}.%
\footnote{A remarkable conclusion from this undoubtedly correct observation
is drawn in Borges'
`argumentum ornithologicum' \citep[p.~29]{Borges1964}.}

Intuitiveness of an operation is explained in
terms of intuitiveness of a proposition:
An operation will be intuitive
if we have intuitive knowledge
\emph{that}
it is well-defined.

As I mentioned,
Parsons argues against the intuitiveness of
the general principle of induction.
His objection turns on its essentially higher-order
character 
\citep[p.~227]{Parsons1986}:
the principle does not speak directly about objects and
operations on them that
are intuitive in his sense,
but about all predicates defined on these objects.
On Parsons' conception,
intuition
\emph{that} is
only possible where all types involved can be instantiated
in (concrete) perception 
\citep[§28]{Parsons2008}.

For the same reason,
Parsons argues that the principle
that 
\emph{every} 
operation on numbers can be iterated
is not intuitive.
The question of its intuitiveness is fact
prior to that for induction,
because,
as Parsons points out,
induction is specific to reasoning
about domains in which the objects are obtained from an
initial object
by arbitrary finite iteration of a given operation.
Similarly,
a simple way of defining an operation inductively
is to define it as the iteration of another operation.
For example,
addition can be defined as iteration of the successor
operation,
multiplication as iteration of addition,
and exponentiation as iteration of multiplication.

This does of course not exclude the possibility
that certain specific operations that are usually
defined inductively are intuitive after all;
it is just that,
in such a case,
this intuitiveness would not have its ground
in the fact that
one way of defining
that operation
proceeds
by appropriately instantiating the
general schema of definition by induction.

An example is addition.
Addition is usually understood according to
an inductive definition like this:
\begin{eqnarray*}
a+0 & = & a\\
a+Sb & = & S(a+b)
\end{eqnarray*}
Parsons accepts an argument proposed by Bernays
that here an alternative understanding
is available
that does not depend on iteration 
\citep[p.~255]{Parsons2008}.
If $a$ is given in intuition (as a string of strokes) and
also $b$,
then so can $a+b$,
because the strings representing $a$ and $b$
can be concatenated in one step,
without iterating through them.

If one abstracts from questions of feasibility,
as Bernays, Parsons, and Brouwer
all do,
it is,
I think,
unproblematic that we can juxtapose
any two strings in one step.
But I am not convinced that this means that
we can come to accept addition
as an intuitive operation
without any appeal to iteration.
We have to verify that the new object
is not a mere juxtaposition of the two original ones,
but itself a string of strokes – a concept
that,
Parsons says,
`involves iteration' 
\citep[p.~175]{Parsons2008}.
In other words,
we have to know that the operation
of concatenating strings is type-preserving.
I do not think that we can know
that the result of an act of concatenation
can in turn be
iterated through
without
having verified this in some particular cases.
(Such cases then serve
to found,
in acts of what Husserl would call
`eidetic variation',
the general judgement.)
But I will leave this aside for the moment.

For multiplication one may propose
an argument analogous to that for addition.
Multiplication is usually understood in a way that depends
on iteration:
\begin{eqnarray*}
a \times 0 & = & 0 \\
a\times Sb & = & a \times b + a
\end{eqnarray*}
However,
given strings $a$ and $b$,
we may obtain $a\times b$ directly,
by replacing each stroke in $b$
by a copy of $a$.
In fact,
we are then doing the same thing as when we
constructed $b$,
except that we now take $a$ as the unit instead of a single
stroke.

Parsons notes
that this way of understanding
addition and multiplication
does not easily generalise
to exponentiation
\citep[p.~256]{Parsons2008}.
The way I would elaborate this is to say that,
in the case of addition and multiplication,
if we leave aside the objection mentioned above concerning
type-preservation,
we can indicate how to transform an image of the function
arguments into an image of the function value
without invoking any arithmetical concept,
let alone iteration of arithmetical operations.
Only direct manipulations of the image strings are
required,
such as copying,
concatenating,
and replacing one string by another.
It seems that all such direct manipulations of strings
can be understood in terms of part-whole relations.
But the relation that exponentiation determines between
one of its arguments,
the exponent string,
and the string representing the result of the operation,
namely,
the former indicates
the number of iterations of multiplication needed to arrive
at the latter,
is not a part-whole relation.

For example,
consider the exponentiation $2^3$.
The salient relation between 3 and 8
that is established by performing this operation
is not that in a construction of
8 strokes,
a string of 3 strokes enters as a
part.
Of course there is a part-whole relation between a
string of 8 strokes and any 3 consecutive strokes in it;
but that relation is not brought about by carrying out
the operation
of exponentiation.
This seems to rule out that
exponentiation can be understood
in the same concrete sense as addition and multiplication
can.%
\footnote{This is how I would begin answering Tait's
criticism  of Parsons' account
(in the earlier presentation of \citet{Parsons1998}):
\begin{quote}
It isn't that we are told [by Parsons] exactly
how to reason about the dead symbols and what the
limits of such reasoning are,
and that what we are told does not support
exponentiation.
Aside from examples,
we are told nothing about how to reason concerning them,
except that it should be logic-free.
It is true that a certain sketch is given of how
to understand addition and multiplication
finitistically,
using the operations of concatenation
and replacement 
(of each occurence of $|$ in a word by a word),
where there is no reasonable extension of this sketch
to exponentiation.
But these constructions are themselves just examples. 
\citep[101n.14]{Tait2010}
\end{quote}}

Bernays
once proposed an argument to the effect that,
even if we accept that any understanding of exponentiation
has to involve an appeal to iteration,
it can still be understood as intuitive on
different grounds.
His argument,
as translated by Parsons,
runs as follows:
\begin{quotation}
\label{BernaysExp}
Consider the example of the number $10^{10^{10}}$.%
\footnote{[Bernays writes $10^{10^{1000}}$.]}
We can reach this number in a finitary way as follows:
We start from the number 10,
which,
in accordance with
one of the normalizations given earlier we represent
by the expression
\begin{center}
$1111111111$.
\end{center}
Now let $z$ be any number which is represented by a
corresponding
expression.
If we replace in the above expression each 1 by the
expression $z$,
then again a number expression arises,
as we can make clear intuitively,
which for communication we designate with “$10 \times z$”.%
\footnote{[Footnote Bernays, omitted by Parsons:
`Here we have a symbol “with content”\,'.]}
Thus we obtain the process of multiplying a number by 10.
From this we obtain the process of passing from a number $a$
to the number $10^a$,
in that we let the number 10 correspond
to the first 1 in $a$ and,
to each attached 1,
the process of multiplying by 10,
and continue until we are at the end of the expression $a$.
The number obtained by the final process of multiplying by
10
is designated by $10^a$.

This procedure offers no difficulty from the intuitive point
of view. \citep[p.~258]{Parsons2008}%
\footnote{\begin{quotation}
\german{Auch hier bestehen Grenzen für die Vollziehbarkeit der Wiederholungen
sowohl im Sinne der wirklichen Vorstellbarkeit wie auch im Sinne
der physikalischen Realisierung.
Betrachten wir beispielsweise die Zahl $10^{10^{1000}}$. Zu dieser können
wir auf finitem Wege folgendermaßen gelangen: Wir gehen aus von der Zahl 10, die wir gemäß der einen von unsern früher angegebenen Normierungen durch die Figur 
\begin{center}
$1111111111$.
\end{center}
repräsentieren. 
Sei nun $z$ irgendeine Zahl, die durch eine entsprechende Figur repräsentiert wird. Ersetzen wir in der vorigen Figur jede 1 durch die Figur $z$, so entsteht, wie wir uns anschaulich klarmachen können, wieder eine Zahlfigur, die zur Mitteilung mit
“$10 \times z$” bezeichnet wird. [Footnote Bernays: `Hier handelt es sich um ein Zeichen 
“mit Bedeutung”\,'] Wir erhalten so den Prozeß der Verzehnfachung einer
Zahl.
Aus diesem gewinnen wir den Prozeß des Überganges von einer Zahl $a$
zu $10^a$,
indem wir der ersten 1 in $a$ die Zahl 10 und jeweils jeder
angehängten 1 den Prozeß der Verzehnfachung entsprechen lassen und hierin so weit gehen,  bis wir mit der Figur $a$ am Ende sind.
Die durch den letzten Prozeß der Verzehnfachung gewonnene Zahl bezeichnen wir
mit $10^a$.}

\german{Dies Verfahren bietet für die anschauliche Einstellung grundsätzlich
 keinerlei Schwierigkeit.}
\citep[pp.~38–39]{Bernays1976}
\end{quotation}}
\end{quotation}
Parsons objects that
`processes are not clearly objects of intuition' \citep[p.~258]{Parsons2008};
I would say that,
from a Husserlian point of view,
they are,
and that,
phenomenologically,
Bernays' account will turn out to be acceptable.
I will come back to the intuitive givenness of processes
below,
and argue that
a phenomenological account
is available
for induction more generally.
All of this obviously requires a notion of intuition that
goes
beyond
the immediate givenness of spatio-temporal configurations,
however idealised.
That is a notion that does not fit
with Hilbert and Bernays' descriptions of intuition
elsewhere;
and not with Parsons' own notion either.

Parsons suggests that Bernays possibly means that
`what the 1
is to be replaced by is
\emph{the result}
of multiplying by
10 what one obtained
at the previous stage'.
Note that this reading of Bernays replaces his appeal to
intuition of a process by an appeal to intuition of
an object (of a type that Bernays accepts as intuitive).
Parsons points out that,
if it is really intuition of the result that is meant,
the argument becomes circular.
It then implicitly appeals to induction:
the argument,
on this reading,
presupposes that,
for any $n$,
when we set out to form
an intuitive representation of $x^{n+1}$,
the result of raising $x$ to the power $n$,
$x^n$,
is intuitively given;
and that assumption is in effect an induction hypothesis.
The problem with this kind of circularity,
Parsons notes,
is not specific to the case of exponentiation;
the problem would arise in any attempt to justify
the intuitiveness of an operation
by defining it as the iteration of another operation
that has already been seen to be intuitive.

\section{Brouwer}
\label{Brouwer}

According to Brouwer,
all mathematical acts
and the objects constructed in them
are developed out of the `basic intuition' or `Urintuition'.
This is based on `the perception of the move of time':
\begin{quote}
the perception of the move of time,
i.e.~of the falling apart of a life moment into two
distinct things,
one of which gives way to the other,
but is retained by memory.%
\footnote{Like Brouwer here,
Husserl at first used the term `memory' 
– specifically, `primary memory' 
(\german{primäre Erinnerung}) – 
in his descriptions of this phenomenon.
However,
as Husserl explains in a key text dated `not before September 1909',
a `memory' always relates to a temporal object already constituted,
whereas what is meant here is rather
an intentional relation between phases
of the flow of consciousness,
which themselves are not temporal objects. 
He therefore introduces a different term for the latter,
`retention' (\citep[p.~333]{Husserl1969}; 
\citep[p.~191]{Husserl1985b} for the date).
That Brouwer in this quotation uses `memory' 
in the sense of retention
is clear from the fact that he invokes it in the context of
the basic perception of the move of time,
which moreover
precedes also in Brouwer's view the appearance of
objects of any kind \citep[p.~1235]{Brouwer1949C}.
Note that Brouwer,
who had conversations with Husserl in 1928,
in his publications never came to adopt the
term `retention' (see, e.g., \citet{Brouwer1952B,Brouwer1954A}).}
If the two-ity thus born is divested of all quality,
there remains the empty form of the common substratum of all
two-ities.
It is this common substratum,
this empty form,
which is the basic intuition of mathematics. 
\citep[p.~141]{Brouwer1952B}%
\footnote{I have chosen this late formulation because it is
well-known and concise.
In his dissertation,
Brouwer formulates the same idea
\citep[p.~8]{Brouwer1907},
as well as in his 1912 lecture
`Intuitionism and formalism'
\citep[p.~12]{Brouwer1912A}/\citep[p.~127]{Brouwer1975}.}
\end{quote}
The `two distinct things' that Brouwer speaks of
are two phases of consciousness,
that of the present and that of the immediate past,
each with their full experiential content.%
\footnote{To apprehend a phase as a thing
in fact requires an act of abstraction
as the phases of consciousness form a continuum,
which is not composed of discrete objects.}
Brouwer elsewhere specifies that the intuitive temporal continuum is
`a measureless one-dimensional continuum in a single subject'
\citep[p.~116]{Brouwer1975},%
\footnote{\dutch{`een maatloos eendimensionaal continuüm in een enkel subject'} \citep[p.~15]{Brouwer1909A}}
the experience of which exists independently 
of any outer experience \citep[pp.~97-98, 118]{Brouwer1907}.
Brouwer's intuitive time corresponds to Husserl's
\emph{inner} time awareness;
Brouwer distinguishes it
from `scientific time',
which is not a priori but a posteriori and
presupposes the existence of mathematics 
developed on the basis of intuitive time 
\citep[pp.~99n.]{Brouwer1907}.
He also emphasises that mathematics,
conceived of as the activity of 
constructing mathematical objects
on its basis of this intuition,
is  not of a linguistic nature 
\citep[pp.~169,176–177]{Brouwer1907}.
It is,
to use his well-known later formulation,
`essentially languageless' \citep[p.~141]{Brouwer1952B}.

For Brouwer,
then,
the objects of
mathematical intuition are not,
as in Parsons' model,
strings of strokes,
but constructions out of
inner time awareness.
But the objects that are intuitive in Parsons' sense can be
mapped
to objects that are intuitive in Brouwer's sense,
as successive strokes in strings can be mapped to
successive intervals in time.%
\footnote{Compare `It is
of course natural
also to view the generation of strings temporally.
I believe that the structure that results,
and the issues concerning it,
are the same as in Brouwer's case.'
\citep[p.~174n.70]{Parsons2008}}
Intuitiveness of operations and of propositions are,
in the abstract,
defined in the same way as for Parsons.
I claim,
but will not argue for it at this point,
that all arithmetical principles that are intuitive
knowledge for
Parsons
are also intuitive knowledge for Brouwer.
On the other hand,
if a certain object or principle is not intuitive on Parsons' account,
it may still be on Brouwer's.

Induction is a case in point.
When,
at the beginning of his dissertation,
Brouwer
gives construction methods
for 
$9\times 4$
and
$4^5$,
these are,
in effect, 
instantiations of
straightforward inductive definitions:%
\footnote{The construction method he gives
for 
$3+4$
is rather of the type
of direct concatenation.}
\begin{quote}
By
$9\times 4$
I mean:
Count up to 4,
write 1 on another line,
add 4 on the first line
(the operation `+4' described above),
write 2 on the second line,
etc.,
till 9 has been written on the second line.
By
$9\times 4$
I mean the last number on the first line.
\citep[pp.~4–5]{Brouwer1907}/\citep[p.~15]{Brouwer1975}
\end{quote}
These are clearly meant to be examples
of constructions that are intuitive,
as is confirmed at the end of that chapter
\citep[p.~77]{Brouwer1907}/\citep[pp.~5–52]{Brouwer1975}.
Heyting,
in his edition of Brouwer's \emph{Collected Works},
comments on this passage that
\begin{quote}
The definitions of the arithmetical operations
by recursion and the derivation of their
properties by induction are intuitionistically correct.
Probably Brouwer intended to demonstrate
that the notions of these operations are more
primitive than the general notions of recursion
and induction.
\citep[p.~565]{Brouwer1975}
\end{quote}
Indeed,
as we will see below,
for Brouwer the general notion of induction and recursion
belong to `mathematics of the second order',
whereas concrete additions and multiplications
belong to `mathematics of the first order'.
But one never finds in Brouwer attempts at accounts
for the intuitiveness of these operations
of the kind proposed by Bernays or Parsons,
and he did accept induction as an intuitive principle.
In the list of `theses' that go with his
dissertation
(\emph{\dutch{stellingen}}),
Brouwer states:
\begin{quote}
The admissibility of complete induction cannot only not be
proved,
but it ought neither to be considered as a separate axiom
nor as a
separately seen intuitive truth.
Complete induction is an act of mathematical constructing,
which is already justified by the basic intuition of
mathematics.
\citep[p.~98, thesis II, trl.~modified]{Brouwer1975}%
\footnote{\dutch{De geoorloofdheid der volledige inductie
kan niet alleen niet worden bewezen,
maar behoort ook geen plaats als afzonderlijk axioma
of afzonderlijk ingeziene intuïtieve waarheid in te nemen.
Volledige inductie is een daad van
wiskundig bouwen,
die in de oer-intuïtie der wiskunde reeds haar
rechtvaardiging heeft.
\citep[loose leaf, \dutch{Stelling} II]{Brouwer1907}/\citep[p.~139]{Dalen2001a}.}}
\end{quote}

Induction ought not to be considered as a separate axiom
because,
in Brouwer's view, 
it 
is neither an axiom nor separate.
It is not an axiom but an act,
and a principle only in the derived sense
of being a correct description of that act,
upon reflection on it.
Induction as a principle would perhaps best be
formulated 
as a rule under which mathematical
construction is closed:
If I have obtained a construction for 
$A(0)$ and
if,
by whatever mathematical 
(not necessarily merely logical)
means, 
I can obtain a construction for
$A(Sn)$
whenever I have obtained a construction for
$A(n)$,
then I have the mathematical means to obtain 
a construction for $\forall xA(x)$.%
\footnote{When in Notebook 1 Brouwer quotes the following passage from
\emph{\french{La science et l'hypothèse}}:
\begin{quote}
\french{Poincaré:
Les mathématiciens procèdent donc « par construction »,
ils « construisent » des combinaisons de plus en plus
compliquées.
Revenant ensuite par l'analyse de ces combinaisons,
de ces ensembles,
pour ainsi dire,
à leurs éléments primitifs,
ils aperçoivent les rapports
de ces éléments et en déduisent les rapports des
ensembles eux-mêmes'}
\citep[p.~26]{Poincare1902}.
\end{quote}
he adds in the margin,
next to `\french{et en déduisent}',
`\dutch{liever}:
\french{et essayent d'en construire}' \citep[Archive, notebook I, p.~36]{Brouwer.Archive}.}
Correspondingly,
an account of induction should be given in terms
not of propositions and operations on them,
but of acts and operations on them.
This of course requires objectification
of acts (see below),
but objectified acts are still different from propositions.

And neither is induction `separate',
because,
as I gloss that term,
it is not evidentially independent.
It is
not an option to do any intuitionistic mathematics
without appealing to something from which,
Brouwer claims,
induction 
can be made evident as well 
– the `basic intuition of mathematics'.
One may of course  restrict oneself to doing
intuitionistic mathematics
without actually using induction,
but there is no analogy between
doing so and,
for example,
doing classical set theory without using the axiom of choice.

A rephrasing of the quoted `thesis' on induction is
known from Brouwer's letter
to Jan de Vries,%
\footnote{1858–1940; professor of geometry in Utrecht from 1897 to 1928. 
No further contact between Brouwer and De Vries before or in the period of Brouwer's 
dissertation is known.
It is likely that the occasion for the exchange had been created by 
Brouwer's thesis advisor Diederik Korteweg,
who was a friend of De Vries and incidentally had been
the thesis advisor of his younger brother Gustav,
who received his PhD in 1894.
See 
\citet{Willink2006}.}
of which the copy 
that remains in the Brouwer archive is undated,
but which was clearly written around the time
of the thesis defence.
The letter comments on various ideas of the
dissertation,
and remarks on induction:
\begin{quote}
I replace the
`axiom of complete induction'
with the
`mathematical construction-act of complete induction'
and show how,
given the intuition of time,
this is nothing new. 
\citep[p.~155]{Dalen2001a}%
\footnote{Brouwer to J.~de Vries: `\dutch{Ik stel in plaats
van het
“axioma van de volledige
inductie” de “wiskundige opbouw-handeling van volledige
inductie”,
en laat zien,
hoe die na de tijdsintuïtie niets nieuws meer is.}'}
\end{quote}

But in spite of what he says here,
no detailed attempt at showing this is found
in the dissertation.
One does however find
an important indication of the form that such an account
should take,
in a footnote to this list of three examples 
of valid (general) synthetic a priori judgements:
\begin{quote}  
\begin{enumerate}
\item the very possibility of mathematical synthesis,
of thinking many-one-ness, and of the repetition thereof
in a new many-one-ness.
\item the possibility of intercalation 
(namely that one can consider as a new element
not only the totality of two already compounded,
but also that which binds them:
that which is not the totality and not the element)
\item the possibility of infinite continuation 
(axiom of complete induction%
\footnote{[Perhaps in a hurry, 
Brouwer here writes 
`axiom of complete induction' 
instead of
`mathematical construction-act of complete induction',
which in the letter to De Vries he says is what he replaces
the former with.]}%
)
\end{enumerate} \citep[p.~70]{Brouwer1975}%
\footnote{\dutch{`1.~de mogelijkheid zelf van wiskundige synthese,
van het denken van veeleenigheid,
en van de herhaling daarvan in een nieuwe
veeleenigheid. 2.~de mogelijkheid van tusschenvoeging, 
(dat men n.l.~als nieuw element kan zien niet alleen
het geheel van twee reeds samengestelde,
maar ook het bindende:
dat wat niet het geheel is,
en niet element is). 3.~de oneindige voortzetbaarheid
(axioma van volledige inductie)'} 
\citep[pp.~119–120]{Brouwer1907}}
\end{quote}
The footnote to this list states:
\begin{quote}
One must however not try to base mathematics or experience
on such judgements:
they are the result of
\emph{viewing} the basic intuition
\emph{mathematically},
and hence presuppose the basic intuition in the viewing as
well in what is viewed;
they belong to what we shall call
in the next chapter
\emph{mathematics of the second order}
\citep[p.~70n, trl.~modified, original emphasis]{Brouwer1975}.%
\footnote{\dutch{Men trachte echter niet,
die oordeelen aan de wiskunde
of aan de ervaring ten grondslag te leggen;
ze zijn het gevolg van 
\emph{wiskundig bekijken}
der oer-intuïtie,
vooronderstellen dus de oer-intuïtie
zoowel in het bekijken als het bekekene;
ze behooren tot wat we
in het volgende hoofdstuk zullen noemen
\emph{wiskunde der
tweede orde}.}
\citep[p.~119n, original emphasis]{Brouwer1907}}
\end{quote}
Mathematics of the second order had,
in fact,
already been defined a few pages earlier:
\begin{quote}
Strictly speaking the construction of intuitive mathematics 
in itself is an action and not a science;
it only becomes a science,
i.e.~a totality of causal sequences,
repeatable in time,
in a mathematics of the second order,
which consists of
the mathematical consideration of mathematics
or of the language of mathematics.%
\footnote{\dutch{Eigenlijk is het gebouw der intuitieve wiskunde zonder meer
een daad,
en geen wetenschap;
een wetenschap,
d.w.z.~een samenvatting van in den tijd herhaalbare causale volgreeksen,
zordt zij eerst in de wiskunde der tweede orde,
die het wiskundig bekijken van de wiskunde of van de taal der wiskunde is.}
\citep[p.~98n1]{Brouwer1907}}
\citep[p.~61n1]{Brouwer1975}
\end{quote}
The latter type of
mathematics of the second order
is the better known one in the literature,
no doubt because it is central to the intuitionistic
conception of logic as the study of patterns
in such descriptions.
But 
in order to make the general principles
of iterability and induction
evident,
it will have to be
the first type that we engage in,
because in this case we are concerned with patterns
in intuitive acts,
not in language.
Although Brouwer points out that
a judgement of second order-mathematics can
play no rôle in founding mathematics,
second order-mathematics
is all the same intuitive,
so that judgements based on it
still express intuitive knowledge.

In his 1911 review of Mannoury's
book \emph{\german{Methodologisches und philosophisches
zur Elementar-Mathematik}},
Brouwer speaks of
`the intuition of complete induction' 
\citep[p.~200]{Brouwer1911A};%
\footnote{\dutch{`de intuïtie der volledige inductie'}}
and
in 1912,
in his inaugural lecture `Intuitionism and formalism',
he remarks that
for finite numbers as understood intuitionistically,
induction is
`evident on the basis of their construction' 
\citep[p.~129–130, trl.~modified]{Brouwer1975}.%
\footnote{`\dutch{dit principe [i.e., volledige inductie],
dat voor de eindige getallen
van de intuïtionist op grond hunner constructie evident is
…}' 
\citep[pp.~15–16]{Brouwer1912A}.}
But this is in passing,
without development either there or in later publications.

In the notebooks in which Brouwer drafted his
dissertation between 1904 and 1907,
induction is occasionally commented on,
but never to question its status of
an acceptable principle.%
\footnote{These notebooks have
neither been translated nor published yet. 
But there are many
quotations from it 
(with translations) in 
John Kuiper's dissertation
\citep{Kuiper2004}.}
One comment,
in the last notebook, 
is of particular interest.
It is occasioned by the
following passage
in Poincaré's \emph{\french{La Science et l'Hypothèse}}:
\begin{quote}
\french{Ce procédé est la démonstration par récurrence.
On établit d'abord un théorème pour $n=1$;
on montre ensuite
que s'il est vrai de $n-1$,
il est vrai de $n$ et on en conclut
qu'il est vrai pour tous les nombres entiers.}
\citep[p.~19]{Poincare1902}
\end{quote}
Brouwer comments:
\begin{quote}
The principle of induction is not:
`if the theorem is correct for 1,
and for \mbox{$n+1$} if it is correct for $n$,
then it is correct for every number',
but the possibility to think the same thing
repeated forever,
so also buildings [i.e.,
construction acts],
so also attempts at buildings
in which at each number one
gets smacked in the face by the 
principle of contradiction. 
\citep[Archive, notebook VIII, p.~65]{Brouwer.Archive}%
\footnote{`\dutch{Het principe der inductie is niet,
dat:
`als stell[ing] geldt voor 1 en voor $n+1$ als voor $n$,
dan voor elk getal',
maar de mogelijkheid,
om zich een gelijksoortig ding altijd door herhaald te
denken,
dus ook gebouwen,
dus ook een poging om te bouwen,
die bij elk getal opnieuw zijn neus stoot door} 
\german{den
Satz vom
Widerspruch}.'}
\end{quote}
It would have been more accurate
for Brouwer to write here,
as he does in the dissertation, 
that 
one is `smacked in the face'
not primarily by the 
(propositional)
principle of contradiction,
but by the fact that a certain construction act 
`does not go through' (`\dutch{niet verder gaat}' 
\citep[p.~127]{Brouwer1907}).
Be that as it may,
what is noteworthy in this quotation is the insistence
on induction as an instance of iteration,
which  Brouwer does not make explicit elsewhere.

Brouwer's preceding remarks on induction,
which are those I have been able to find,
may be summarised in
five tenets:
\begin{enumerate}
\label{B1B5}
\item[B1] induction is primarily an act, not a proposition;
\item[B2] the act of induction is an instance of iteration;
\item[B3] the intuition of time is a condition of possibility of the act of induction;
\item[B4] the judgement of the validity of the induction principle is a result of `second-order mathematics';
\item[B5] the induction principle is evident on the basis of the intuitionistic construction of the finite numbers.
\end{enumerate}
Before turning to the details of an account based on these tenets,
which is not to be found in Brouwer's published
or unpublished writings,
I should like to show to what (varying) degree
several published intuitionistic accounts of later authors
do not show these tenets B1–B5.

I do not know how much of
(the ideas in)
Brouwer's unpublished notes
on induction were known to his foremost student,
Arend Heyting.
But Brouwer's emphasis on acts and iteration
over
propositions and deduction
is absent from the 
explicit justification of induction
proposed by Heyting in his
\emph{Intuitionism. An Introduction} \citep{Heyting1956}.
Brouwer's claim in the inaugural lecture that induction is
`evident on the basis of the intuitionistic construction of the finite numbers',
however, 
is echoed there,
although one suspects that Brouwer and Heyting had different ideas
as to exactly how induction becomes evident on that basis:%
\footnote{Heyting had published the same justification in Dutch in 1936,
in his paper \citet{Heyting1936}.
In the main text I quote from the book of 1956 because it is much better known.}
\begin{quote}
Clearly the construction of a natural number $n$
consists in building up successively certain natural
numbers,
called the numbers from 1 to $n$,
in signs:
$1 \rightarrow p$.
At any step in the construction we can pause
to investigate whether the number reached
at that step possesses a certain property or not.
…%
\footnote{The omitted part is the example of
an argument for the theorem that $m \neq n$
and $m>n$ implies $n<m$.}
The theorem of complete induction admits of a proof
of the same kind. 
Suppose $E(x)$ is a predicate of natural numbers
such that $E(1)$ is true and that,
for every natural number $n$,
$E(n)$
implies $E(n')$,
where $n'$
is the successor of $n$.
Let $p$ be any natural number.
Running over $1 \rightarrow p$ 
[the number 1 to $p$ in their natural order]
we know that the property $E$,
which is true for 1,
will be preserved at every step in the construction of $p$;
therefore $E(p)$ holds.%
\end{quote}
Justifications of the same type are given
in Troelstra's
\emph{Principles of Intuitionism}
of 1969 
\citep[p.~12]{Troelstra1969}
and in Troelstra and Van Dalen's
\emph{Constructivism in Mathematics}
of 1988.
To quote the latter:
\begin{quote}
The justification of induction is based on the mental picture
we have of the natural numbers,
obtained by successively adding `abstract units';
given $A(0)$,
$\forall x(A(x)\rightarrow A(Sx))$
we build parallel to the construction of
$n \in N$ a proof of $A(n)$:
\begin{prooftree}
\AxiomC{$A(0)$}
\AxiomC{$A(0)\rightarrow A(1)$}
\BinaryInfC{$A(1)$}
\AxiomC{$A(1)\rightarrow A(2)$} 
\BinaryInfC{\parbox{4em}{$A(2)\ \dots$\\ \hphantom{$A(2)$}etc.}}
\end{prooftree}
\citep[vol.1, p.~114]{Troelstra.Dalen1988}
\end{quote}
While this does suffice to show
that
starting from the construction of a given number $n$,
we can obtain a proof of $A(n)$,
it does not suffice for demonstrating that
we have one construction that works for all 
(infinitely many)
$n$ at once,
as would be required for a justification of complete induction.
As Parsons notes in his discussion
of this type of reasoning (without citing a particular occurrence),
to claim that it does yield complete induction is fallacious:
\begin{quote}
In fact,
for each $x$,
we can construct a formal proof of $A(x)$ by beginning with
$A(0)$
and building up by modus ponens,
using $A(x) \rightarrow A(Sx)$.
As a
\emph{proof}
of induction, this is circular:
the `construction' of $x$ by a
\emph{succession}
of steps is itself inductively defined,
and it is by a corresponding induction
that it is established that $A$ holds
at each point in the construction. 
\citep[p.~266, original emphasis]{Parsons2008}%
\footnote{Similarly,
Yessenin-Volpin 
rejected both mathematical induction and the soundness principle that
`If the axioms of a formal system are true and the rules of inference
conserve the truth then each theorem is true',
for the reason that these are mutually dependent.
`It is essentially on these grounds that I am not searching for any
axiomatic theory in my program.' \citep[p.~6]{Yessenin-Volpin1968}.
Note that Brouwer's grounding of induction,
discussed below,
is not axiomatic.}
\end{quote}
Dummett argues
in
\emph{Elements of Intuitionism} 
that,
although there is `no uniform proof skeleton' 
for proofs of $A(n)$ from the induction basis
and the induction step,
except if one presupposes induction,
we recognise all the same that
the operation of chaining,
at each $n$,
$n$ applications of modus ponens,
will yield a proof of $A(n)$ for each $n$
\citep[p.~9]{Dummett2000b}.
I do not see how that argument fares better.

The circularity that arises can
also be analysed in terms of the BHK explanation,
which construes an implication as the existence of a
function from proof objects of 
the antecedent to proof objects of the consequent.
For each $n$,
the instantiation
of the premise
$\forall x (Ax \rightarrow ASx)$ required for
the application of modus ponens
yields a different function
(one with as domain proofs of $A(n)$
and as range proofs of $A(Sn)$).%
\footnote{In the proof tree given by Troelstra
and Van Dalen, 
shown above,
these steps of instantiation are left implicit.} 
These infinitely many different functions 
cannot be used directly in a finite proof;
one has to use induction to reason about the application of all of them,
and this introduces the circularity.

As we saw,
Brouwer's idea of the relation
between the evidence for
induction and the construction of the natural numbers
was different from its construal in this circular argument;
for him the central idea was iteration.
In 1960 a constructive account was proposed that did not,
as its originator might have put it,
(cl)aim to reconstruct Brouwer's thought,
yet in effect agreed with
Brouwer
on this point
\citep{Kreisel1960}.
It was part of Kreisel's `theory of constructions',
and runs as follows 
($\star$ is the concatenation operator on constructions):
\begin{quotation}
\french{Supposons qu'une propriété $P$
(portant sur les nombres naturels)
soit déterminée par la construction
$\rhop(b,a)$
($=0$ si $b$ est une preuve que 
$a$ 
satisfait $P$,
$=1$ dans le cas contraire)
et que les résultats
\[
P(0),
P(a)\rightarrow P(a \star  1)
\]
soient établis;
autrement dit,
on a deux constructions $b_0$ et $\rho$
telles que
$\rhop(b_0,0)=0$
et que
\{%
$\rhop(b,a)=0$
implique
$\rhop[\rho(b,a),a \star  1]=0$\}.}

\french{Alors,
étant donnée une construction 
$a$
faite à partir de (i) et (ii),%
\footnote{[Defined on p.~389 as (i) $Z(0)$ and (ii) $Z(a)\rightarrow Z(a\star 1)$ (which
formulations incidentally go back to Hilbert)]}
il suffit de
\emph{suivre}
cette construction pour obtenir une
$b_a$ 
telle que 
$\rhop(b_a,a)=0$:
à chaque application de (ii)
dans la construction de 
$a$
correspond une application de
$\rho$.}

\french{Ce raisonnement s'exprime dans le cadre de
la théorie des constructions abstraites,
qui a pour base certains axiomes existentiels
assez élémentaires.}
\citep[p.~390]{Kreisel1960}%
\footnote{Also, in a slightly different form, in
\citet[p.~209]{Kreisel1962}.}
\end{quotation}
This account is of a different type than Heyting's.%
\footnote{In the 1962 paper mentioned to in the previous note,
Kreisel does refer to Heyting 1956,
for its informal explanation of the meaning of the intuitionistic
logical constants 
\citep[p.~98]{Heyting1956};
but Kreisel does not comment on Heyting's argument for induction 
elsewhere in that book.}
It in effect avoids circularity by 
appealing to one operation,
$\rho$,
instead of infinitely many.
It achieves this by taking
not just a natural number as argument,
but also an abstract construction, 
in such a way that the operation can be iterated.
These `abstract constructions' are,
in Kreisel's framework,
themselves objects of which only their formal properties are
taken into account;
it is 
a `formal semantic foundation'
the value 
of which,
qua formal theory,
Kreisel acknowledges to be 
`primarily technical'
\citep[p.~199]{Kreisel1962}.
The questions
if and how his notion of abstract construction
may be related to Brouwer's notion
of construction objects as resulting
from construction acts out of a `basic intuition'
thus fall out of the scope of Kreisel's paper.%
\footnote{Kreisel acknowledges that he indulges in 
`the “mixing” of mathematics and metamathematics stressed
in the informal writings of intuitionists' 
\citep[p.~202n.9]{Kreisel1962}
– stressed,
one may add,
as something that goes against the foundational order
\citep[pp.~169–178]{Brouwer1907}.}

An approach that likewise leads to a construal of induction as iteration,
but now in explicitly Brouwerian terms and 
(therefore)
not in the context of
formal semantics,
was given by Van Dalen in 2008:
\begin{quote}
Given $A(1)$ and \mbox{$\forall n
(A(n)
\rightarrow
A(n+1))$},
we
want to show \mbox{$\forall n (A(n))$}.
`Show' means for a
constructivist `present a proof',
where we have to keep in
mind that
already in 1907 Brouwer was aware that proofs are
constructions;
he spoke of `erecting mathematical buildings' and `fitting
buildings into other buildings'.
In modern terms this would be read as
`constructing mathematical structures' and `constructing a
structure on the basis of
(out of)
another structure'.
It is quite clear that
he knew how proof-constructions for implication,
universally and
existentially quantified statements were to be made.
The cases of
conjunction and disjunction were tacitly understood.
So
– returning to the matter of induction –
we may assume that there is a proof $a_1$ of $A(1)$;
notation – $a_1 : A(1)$. 
Now a proof for \mbox{$\forall n
(A(n) \rightarrow A(n+1))$}
is a construction $c$ that for any given $n$ and proof $a:
A(n)$
yields a proof \mbox{$c(n,a) : A(n+1)$}.
So $a_2 = c(1,a_1)$ and $a_2 : A(2)$, and $a_3 = c(2,a_2)$
and $a_3:A(3)$,
… 
Hence parallel to the construction of the
natural numbers we obtain the
(potentially infinite)
sequence of proofs
$a_1, a_2, a_3, \dots $, i.e., a proof for $\forall n (A(n))$.
\citep[p.~11]{Dalen2008}
\end{quote}
In a footnote,
Van Dalen adds that
`In systems with an explicit recursor, one can often
write down a term for the proof-construction given by the
sequence' \citep[p.~11n.6]{Dalen2008}.
The procedure sketched
may be thought of in terms of 
a primitive recursive function $f$ such that
$f(n)$ is a proof $a_n$ of $A(n)$:
\[
\begin{array}{lll}
f(1) & = & a_1\\
f(Sn) & = & c(n,f(n))
\end{array}
\]
In a type theory with a recursor $\mathrm{R}$
satisfying the conversion relations
\[
\begin{array}{lll}
\mathrm{R}uv1
& \leadsto & u \\
\mathrm{R}uv(Sn) & \leadsto & vn(\mathrm{R}uvn)
\end{array}
\]
we can put $a_1$ for $u$,
and for $v$ an operation corresponding to $c$.
In
(slides for)
a lecture in Groningen in 2009,
Van Dalen claims that
\label{Recursor}`The Ur-intuition also yields the recursor!'.%
\footnote{`\dutch{De herhaalbaarheid van de
successor-operatie is een onderdeel van de oer-intuïtie, de
z.g.~“zelfontvouwing”. De oer-intuïtie levert ook de
recursor!}'
\citep[slide 25]{Dalen2009}
}

I understand that claim as follows.
Brouwer in 1907 of course did not have a theory of recursive functions,%
\footnote{In a letter of July 17, 1928, Brouwer suggested that Heyting
add to the latter's formalisation of intuitionistic logic and analysis
a formalisation of the notion of law 
(in the sense of a so-called spread law)
\citep[p.~334]{Dalen2011}; 
Heyting did
not take this up.
A lawlike sequence,
or to be precise a spread with a lawlike sequence as its single element,
is given as a limiting case;
and a recursive sequence is intuitionistically lawlike.
(Whether the converse is also true is a different question.)}
but he did,
as we have noted,
have a solid idea of iteration and
its relation to induction.
Primitive recursion can be reduced
to iteration;
applied to the recursive function above,
this can be done
as follows.%
\footnote{In the literature, 
the reduction of recursion to iteration goes back to Kleene's
iterative rendering of the recursive predecessor function
\citep{Kleene1936}. For a systematical discussion,
see \citet{Robinson1947}.}
In a type theory with an iterator $\mathrm{I}$
\[
\begin{array}{lll}
\mathrm{I}uv1
& \leadsto & u \\
\mathrm{I}uv(Sn) & \leadsto & v(\mathrm{I}uvn)
\end{array}
\]
and with pairing,
we can set
\begin{equation*}
u=\langle 1,a_1\rangle
\end{equation*}
and
\begin{equation*}
v(\mathrm{I}uvn) =\langle S\pi_1(\mathrm{I}uvn),  c(\pi_1(\mathrm{I}uvn),\pi_2(\mathrm{I}uvn)) \rangle
\end{equation*}
where $\pi_1$ and $\pi_2$ are the left and right projection operators.
Note that the operations of pairing and projection
are readily understood in terms of
an invocation of
Brouwer's two-ity.
Ordered pairing is acceptable as a general intuitive
operation because it consists
in an order-preserving mapping of the two parts of the empty two-ity
onto whatever the elements of the pair will be,%
\footnote{Brouwer sees the act of `taking together' as the mapping
of the two-ity onto the things that are (to be) taken together
\citep[p.~179n1]{Brouwer1907}/\citep[p.~97n1]{Brouwer1907}.}
and projection in separating one element out of a two-ity.

In the following,
I should like to argue in some detail
that
Brouwer,
in holding to the tenets
B1-B5 (p.~\pageref{B1B5} above),
indeed had in mind
a development that leads
from time awareness to intuitiveness
of the iterator,
and, 
from there,
of induction.

Picking up from the two-ity again,
we first turn to the natural numbers.
Brouwer identifies the natural number 2 with the empty two-ity.
Once I have created an empty two-ity as an object,
time moves on again,
a created two-ity sinks into the past and this,
when I decide to turn my attention to it,
will then
become
one component of a new two-ity.
Brouwer identifies this new,
nested
two-ity with the natural number
3.
Because time keeps moves on,
and I can keep turning my attention to it,
I thus obtain an
\emph{iterative}
structure,
and thereby the natural numbers.
There is an intrinsic ordering of the natural numbers as
constructed intuitionistically
in the sense that
the construction object of $n+1$
includes that of $n$ as a proper part.

Brouwer calls the successive construction of
iterated two-ities the `self-unfolding' of the empty
two-ity.
The idea of this self-unfolding includes the ideas that
time is potentially infinite,
and that in particular every life moment will fall apart
into a two-ity,
one part of which is a previous two-ity.
These ideas obviously cannot be taken to be
empirical observations,
but should be seen as a priori insights into the
structure of time awareness.%
\footnote{See also Parsons' remarks on this point   
\citep[pp.~177–178]{Parsons2008}.}

In some of Husserl's writings,
the relation between iteration and the continuum of
time awareness is made explicit.
Central to his analyses of inner time awareness are
the notions of
retentional and protentional intentionality 
– the intentional directedness
of the living present towards the phase of 
the stream of consciousness that has just elapsed
and towards the phase that is about to come,
respectively \citep[p.~297, p.~333]{Husserl1969}.
In a text of 1916,
published as Appendix I 
to the 1905 lectures on time awareness,%
\footnote{As related by John Brough,
Rudolf Bernet dates the original manuscript to 1916
\citep[p.~105n2]{Husserl1991}.}
Husserl wrote:
\begin{quote}
It is inherent in the essence of
every linear continuum that,
starting from any point whatsoever,
we can think of every other point as
continuously produced from it;
and every continuous production is a production
by means of continuous iteration.
We can indeed divide each interval
in infinitum and,
in the case of each division,
think of the later point of the division
as produced mediately through the earlier points …
Now this is also true in the case of temporal
modification – or rather,
while the use of the word `production'
is a metaphor in the case of other continua,
here it is used authentically.
The time-constituting continuum is a flow of continuous production
of modifications of modifications.
The modifications in the sense of iterations
proceed from the acutally present now,
the actual primal impression 
\emph{i}… 
\citep[p.~106]{Husserl1991}%
\footnote{`\german{Im Wesen jedes linearen Kontinuums liegt es,
daß wir,
von einem beliebigen Punkt ausgehend,
jeden anderen Punkt aus ihm stetig erzeugt denken können,
und jede stetige Erzeugung ist eine Erzeugung durch
stetige Iterierung.
Jeden Abstand können wir ja in infinitum teilen und
bei jeder Teilung den späteren Teilungspunkt mittelbar
durch die früheren erzeugt denken … 
So ist es nun auch bei der zeitlichen Modifikation,
oder vielmehr,
während sonst,
bei anderen Kontinuis,
die Rede von der Erzeugung ein Bild ist,
ist sie hier eine eigentliche Rede.
Das zeitkonstituierende Kontinuum ist ein Fluß
stetiger Erzeugung von Modifikationen von Modifikationen.
Vom aktuellen Jetzt aus,
der jeweiligen Urimpression u,
gehen die Modifikationen im Sinn von Iterationen …}' 
\citep[p.~100]{Husserl1969}}
\end{quote}%
Likewise, in a manuscript of 1931,
Husserl speaks of the `the future horizon as 
flowing continuity of the implicit iteration
of coming realisations of always new presents'
and `this flowing continuity of pasts,
an endless horizon of iteratively nested pasts'
\citep[p.~405]{Husserl2006}.%
\footnote{`\german{[der] Zukunftshorizont als strömender Kontinuität der
impliziten Iteration kommender Verwirklichungen von immer
neuen Gegenwarten}' 
and 
`\german{diese strömende Kontinuität der Vergangenheiten,
[ein] endloser Horizont wiederum iterativ ineinander
geschachtelter Vergangenheiten}'.}

Thus, 
the iterative structure of nested two-ities arises
as an abstraction from a structure present in
inner time awareness:
In the perception of the move of time,
just as the present phase of consciousness 
is modified into a retention,
prior retentions are modified into retentions of themselves,
so that,
for example,
a simple retention now becomes a retention of a retention.
In the constitution of a two-ity
out of two phases of consciousness,
the form of the retentional relation between these two phases
is retained,
while abstraction is made from 
the intermediate retentional relations along the temporal continuum in between.
Protention,
in turn,
is an intention directed towards
the phase of consciousness that is just about to come.
It plays a rôle in 
the
intuition that the sequence of numbers constructed so far 
can always be extended.%
\footnote{For a phenomenological constitution analysis of
potentially infinite sequences 
(with choice sequences as a special case),
see \citet[section 6.2]{Atten2007}.}
I read the statement in Brouwer's notebooks that
`The sequence $\omega$ can only be constructed on the
continuous
intuition of time'
\citep[Archive, notebook IX, p.~29]{Brouwer.Archive}%
\footnote{`\dutch{De reeks $\omega$ is alleen op te bouwen op de
continue
tijdsintuïtie.}'}
as,
in effect,
a recognition of this dependence of the
construction of nested two-ities on the
retentional and protentional intentionalities
that characterise the awareness of inner time.

We obtain intuitive knowledge of the structures of time
awareness is through acts of reflection
– reflection
in the phenomenological sense of
turning our attention to an earlier episode in our flow of
consciousness.%
\footnote{Husserl critically discusses various skeptical claims
regarding reflection in section 79 of Ideen I \citep{Husserl1976a}. 
See also \citet{Hopkins1989}.}
The structure of inner time awareness is not given passively
altogether,
for we have to engage in the appropriate kind of reflection,
which is an activity.
But that activity consists in objectifying the
flow of time and its parts and phases;
it does not
\emph{form}
that which is reflected
on.
There is no circularity,
then,
in saying that iterative structures
are constructed in intuition by projecting
from the iterative structure of time awareness,
as we can obtain intuitive knowledge of whatever
structure inner time awareness exhibits
without first having had to engage in an activity of 
constructing it.
Note that this means that the form of inner time
is not categorial form in Husserl's sense,
as he holds that categorial form is constituted in 
\emph{active formation}
by the ego.
Indeed,
although Husserl for some time between
1907 and 1909 did believe that the form of time is
categorial,
he gave up that idea
when he around 1909
discovered the absolute flow of time which constitutes
itself as a flow – one is tempted to say,
which unfolds itself.%
\footnote{For more on this, see \citet{Atten2015b}.}

The question may be raised how we know that
the construction of the natural numbers
based on these ideas of Brouwer and Husserl
results in the standard numbers.
After all,
there are constructive proofs
of the existence of non-standard models of the theory of the
successor
and,
more importantly,
of Heyting Arithmetic (HA).
Dedekind already devised
a non-standard model for the theory of successor
\citep[p.~100]{Dedekind1890}.
As for non-standard models of arithmetic,
the existence proofs by
Skolem and Gödel are of course not constructive,
and McCarty 
has proved constructively that,
under the assumption of Markov's Principle
and a weak form (consequence) of Church's Thesis,
HA has no
non-standard models \citep{McCarty1988}.
The latter result does not settle the matter however
as,
in a Brouwerian setting
including the theory of the creating subject,
even the weak form of Church's Thesis involved
is false,
and there is a weak counterexample to Markov's Principle.%
\footnote{That weak form of Church's Thesis is 
\[
\forall n(P(n) \vee \neg P(n)) \rightarrow 
\neg\neg\exists e\{e\} \textrm{ is $P$'s
characteristic function}
\]
This is inconsistent with Kripke's example of a function that the creating subject
is able to compute yet cannot be assumed to be recursive, on pain of contradiction 
\citep{Atten2008a}.
Markov's Principle is used in this form:
\[
\forall n(P(n) \vee \neg P(n)) \rightarrow 
(\neg\neg\exists n P(n) \rightarrow \exists n P(n))
\]
As pointed out in Troelstra and Van Dalen 
\citep[I:pp.~205–206 and 237]{Troelstra.Dalen1988},
this is equivalent to
\[
\forall x\in \mathds{R}(x \neq 0 \rightarrow \exists k(|x|>2^{-k}))
\]
to which Brouwer presented a weak counterexample in
`Essentially negative properties'
\citep{Brouwer1948A}.}
Indeed,
De Swart's proof of the compactness theorem
for intuitionistic predicate logic 
\citep[section 3]{Swart1977}
does give rise
to non-standard models 
(which have not been studied so far)
acceptable from a Brouwerian point of view.%
\footnote{A constructive non-standard model that,
unlike De Swart's construction,
does
not depend on choice sequences,
and which has given rise to further work,
was later built by  Moerdijk \citeyearpar{Moerdijk1995}.}

However,
non-standard models,
even if constructive, 
pose no threat to
the intuitionistic account of the natural numbers and their arithmetic.
Dummett has written,
in a passage that Parsons 
\citeyearpar[p.~279]{Parsons2008}
draws attention to:
\begin{quote}
Within any framework which makes it possible to speak coherently
about models for a system of number theory,
it will indeed be correct to say that there is just one
standard model,
and many non-standard ones;
but since such a framework within which a model for
the natural numbers can be described will itself
involve either the notion of `natural number'
or some equivalent or stronger notion such as `set',
the notion of a model,
when legitimately used,
cannot serve to explain what it is to know
the meaning of the expression `natural number'.
\citep[p.~193]{Dummett1978}
\end{quote}
Within the specifically intuitionistic context,
the point can be strengthened
by observing that 
the natural numbers are
privileged not only conceptually,
but also genetically.
The construction of any model 
of the theory of successor or of arithmetic,
whether standard or non-standard,
depends,
like all mathematical activity,
on the unlimited self-unfolding
of the empty two-ity,
and hence on the intuition `and so on'.
But if one
acknowledges that intuition,
one obtains,
by thematising a structure that is already present in time
awareness,
the natural numbers directly from it.

Parsons notes that
\begin{quote}
the point of Dummett's observation that the notion of natural number
must be used in the construction of models of arithmetic is that,
in the end,
we have to come down to mathematical language as 
\emph{used},
and this cannot be made to depend on semantic reflection on that same language.
\citep[pp.~287-288]{Parsons2008}
\end{quote}
Phenomenologically,
one would take certain properties of the mind
to be (partly)
explanatory of the constraints on the use of language
that Dummett takes as the point of departure
for his meaning-theoretical considerations.
For example,
that we have the capacity to iterate the construction
of two-ities is an intuitive truth,
as is the fact that,
with our kind of mind,
we cannot complete infinitely many iterations;
Dummett's observation that
\begin{quote}
Even if we can give no formal characterisation
which will definitely exclude all such [non-standard]
elements,
it is evident that there is not in fact any possibility
of anyone's taking any object,
not described (directly or indirectly)
as attainable from 0 by iteration of the successor
operation,
to be a natural number. 
\citep[p.~193]{Dummett1978}
\end{quote}
depends on
these facts.
A vagueness remains in the circumstance
that,
although at any given moment in the
self-unfolding of the two-ity,
in a sense only finitely many
unfoldings have been made,
this sense of `finite' cannot
be replaced by a definition in numerical terms,
as then the
account would become circular.%
\footnote{This is also pointed out by Dedekind, in
section 6 of his letter to Keferstein \citep{Dedekind1890}. 
In effect,
Dedekind himself
relies on such a pre-numerical understanding of `finite' in his understanding
of proofs as finite objects.}
So this particular notion of
finiteness has to be taken as understood prior to the notion
of natural number.

As quoted above,
Brouwer in his dissertation claims that the
a priori judgement that 
mathematical synthesis can be repeated
is a result of `mathematics of the second order',
the result of viewing mathematics mathematically.
In that same work,
Brouwer defines 
mathematical viewing as 
`seeing repetitions of sequences';%
\footnote{\dutch{wiskundig bekijken}, \dutch{het zien van herhalingen van
volgreeksen} \citep[pp.~81, 105-106]{Brouwer1907}.}
in 1929 he analyses it into the two phases of assuming the 
`temporal attitude'
– accepting inner time awareness,
which is a necessary condition for mathematics – 
and then the `causal attitude':%
\footnote{`zeitliche Einstellung'; `kausale Einstellung'}
\begin{quote}
Nunmehr besteht die
\emph{kausale Einstellung}
im Willensakt der `Identifizierung' verschiedener
sich über Vergangenheit und Zukunft erstreckender
zeitlicher Erscheinungsfolgen.
Dabei entsteht ein als
\emph{kausale Folge}
zu bezeichnendes gemeinsames Substrat dieser
identifizierten Folgen.
\citep[p.~153]{Brouwer1929A}
\end{quote}

If we count the givenness of an objectified act in reflection as an 
`\german{Erscheinung}',
then it is clear that viewing our mathematical activity mathematically
allows for the total or partial identification of temporally distinct acts
(or series of acts)
of mathematical construction;
this is a way of `seeing repetition'.
Such identification is a form of applied mathematics,
first,
because taking two (or more) objectified acts together
in one awareness
depends on the formation of two-ities of them,
and second,
because isolating a mathematical structure that is common to two acts
consists in the construction and then successful projection
of the same mathematical structure onto both objectified acts.
As construction and projection themselves take place
in intuition,
the determination of a common structure
leads to intuition of act types 
and,
correspondingly,
of
object types,
namely the type of object constructed in acts of a given type.
Thus,
instead of having to keep,
e.g., 
empty two-ities constructed at different times
ontologically distinct,
we can identify them and consider these as repeated constructions
of the same empty two-ity.
It is also in reflection that the `processes'
in Bernays'
account of exponentiation 
(quoted on p.~\pageref{BernaysExp})
can be given in intuition as individual objects 
and,
founded on that,
as types.

Partial identification is the foundation for abstraction.
For example,
in reflection on acts in which we perform
$5+2$ and $7+2$,
we may come to identify the act type `adding 2 to a natural number'.
This is the example that
came up at Brouwer's
thesis defence,
where the question of the possibility of identifying act types intuitively
was raised.
Gerrit Mannoury there
objected that
in `and so on',  
the `so' 
is not primitive
(and hence neither is `and so on'),
but consists in a relation between relations 
\citep[p.~151]{Dalen2001a}.%
\footnote{See also Brouwer'sremark in a letter to Korteweg of
January 16, 1907, 
on a letter in which Mannoury,
as was the custom,
had informed Brouwer of his
planned objection.
Mannoury's letter unfortunately seems
not to have been preserved.}
Mannoury's idea seems to have been
that
the relation between 5 and 7
and the relation between 7 and 9
are,
in their full intuitive concreteness,
different;
but a similarity relation holds between them. 
Brouwer replied:
\begin{quote}
What you say at the end,
namely that the `so' in `and so on'
is just a relation between relations
and not itself a relation,
can,
it seems to me,
not be upheld either.
Mathematics could not exist,
if I cannot repeatedly think
\emph{the same} 
thing again,
e.g.,
first Jan,
then Piet,
then the same Jan again.
Likewise I can think
\emph{the same}
relation,
and so 
between five and seven
there exists
after all
the same relation as
between seven and nine,
namely, `$+2=$'.%
\footnote{`\dutch{Wat U tenslotte zegt,
dat enzovoort in
\emph{zo}
slechts een relatie tussen relaties 
en niet een relatie zelf ziet,
is geloof ik evenmin vol te houden.
Wiskunde zou niet kunnen bestaan,
als ik niet meermalen weer 
\emph{hetzelfde}
ding kon denken,
b.v.~eerst Jan,
dan Piet,
dan weer diezelfde Jan.
Zo kan ik ook meermalen
\emph{dezelfde}
relatie denken,
en zo bestaat wel degelijk
tussen vijf en zeven
dezelfde relatie als
tussen zeven en negen
nl.~`$+2=$'.}
\citep[p.~151]{Dalen2001a}}
\end{quote}
Unfortunately,
this is where Brouwer's reply ends.%
\footnote{Compare Tait in `Finitism':
`we understand 
$n+2$
not via understanding each of the infinitely many instances,
$0+2$,
$1+2$,
and so on.
Rather,
we understand these via our understanding of what it is
for one sequence to be a two-element extension of another.'
\citep[p.~530]{Tait1981}}
But
this insight naturally leads
from the generation of the natural numbers,
that is,
from the iteration of the successor operation,
to the general principle of iteration.
First observe that,
on
Brouwer's conception of intuition,
it is not just the natural numbers
as individual objects that are constructed in intuition,
but so is the (growing)
object that is the potentially infinite sequence of them.
He speaks of
\begin{quote}
the intuitive truth that
mathematically we cannot construct but
finite sequences,
and also,
on the basis of the clearly conceived `and so on',
the order type ω,
but only consisting of
\emph{equal elements}.
\citep[p.~80, trl.~modified]{Brouwer1975}%
\footnote{`\dutch{de intuitieve waarheid,
dat wij wiskundig niet anders kunnen scheppen,
dan eindige rijen,
verder op grond van het duidelijk gedachte `en zoo voort'
het ordetype ω,
doch alleen bestaande uit
\emph{gelijke elementen}}'
\citep[p.~142–143]{Brouwer1907}.}
\end{quote}
A footnote elucidates the `and so on':
\begin{quote}
\label{L1}The expression `and
\emph{so} on' means the indefinite repetition
of
\emph{one and the same}
object or operation,
even if that object or that operation is defined in
a rather complex way.
\citep[p.~80n1, original emphasis]{Brouwer1975}%
\footnote{\dutch{`Waar men zegt:
“en
\emph{zoo} voort”,
bedoelt men het onbepaald herhalen van
\emph{eenzelfde}
ding of operatie ook al is dat ding of die operatie
tamelijk complex gedefinieerd.'}
\citep[p.~143n1, original emphasis]{Brouwer1907}}
\end{quote}
Now suppose we have a function 
$f:\mathds{N}\rightarrow \mathds{N}$
(or,
more generally,
$f:A\rightarrow B$
such that
$B\subseteq A$).
We now
construct,
in intuition, 
two potentially infinite sequences in parallel,
one being that of the natural numbers starting at 1,
and the other
the sequence
of (the results of) the operation of applying $f$,
beginning with $f(a)$ for some given $a$:
\[
\begin{array}{rcl}
1 & & f(a)\\
2 & & f(f(a))\\
3 & & f(f(f(a)))\\
\vdots & & \vdots
\end{array}
\]
The object constructed in these acts in parallel with $n$ is $f^n(a)$.
The justification for this claim is that,
on the hypothesis that we have a construction method for
a given natural number $n$,
we also know
that the series of
applications of $f$ of length $n$
admits of composition,
because each time we apply one and the same operation
$f$ whose range is included in its domain.
Moreover,
by an appeal to the facts
that the sequence 
construction of natural numbers is given as such
in intuition,
and that the construction of $f^n(a)$ proceeds in
parallel with that
of $n$,
we also know that we have an intuitive construction of
the potentially infinite sequence of
the $n$-fold iterations.

As the insight that operations 
$\mathds{N}\rightarrow \mathds{N}$
can be iterated depends
immediately on the insight that the construction of the
natural numbers is iterative,
the former is hardly more reducible than the latter.
At the end of the dissertation,
Brouwer emphasises:
\begin{quote}
there
are elements of mathematical construction
that
in the system of definitions
which must remain irreducible,
and which therefore,
when communicated,
must be understood from a single word,
sound, or symbol;
they are the elements of construction that are
immediately read off from the Urintuition
or intuition of the continuum;
notions such as
\emph{continuous},
\emph{unity},
\emph{once more},
\emph{and so on}
are irreducible.%
\footnote{\dutch{`Er zijn elementen van wiskundige bouwing,
die om het systeem der definities onherleidbaar moeten
blijven,
dus bij mededeling door een enkel woord,
klank of teken,
weerklank moeten vinden;
het zijn de uit de oer-intuïtie of de
continuum-intuïtie afgelezen bouwelementen;
begrippen als
\emph{continu},
\emph{eenheid},
\emph{nog eens},
\emph{enzovoort}
zijn onherleidbaar.' [1907:180]}}
\end{quote}

Brouwer's view of induction as an act
in mathematical intuition
depends on his particular
understanding of a mathematical property
not in terms of definitions of predicates
but in terms of mathematical construction acts,
which for him are languageless:
\begin{quotation}
\label{BrouwerProperty}
Often it is very simple to introduce inside such a
system,
independently of the way it came into being,
new buildings,
as the elements of which we take elements of the old
one or systems thereof,
in a new arrangement,
but bearing in mind the arrangement in the old building.
What are called the `properties' of a given system
amount to the
possibility of building such new systems in a
certain connection with a priorly given system.

And it is exactly this
\emph{fitting in}
of new systems
\emph{in a given system}
that plays an important part in building
up mathematics,
often in the form of an inquiry into the possibility or
impossibility of a fitting-in satisfying certain conditions,
and in the case of possibility into the various ways in
which it is
possible. 
\citep[p.~51–52, trl.~modified, original emphasis]{Brouwer1975}%
\footnote{\begin{quotation}
\dutch{Binnen zulk een opgebouwd systeem zijn dikwijls,
geheel buiten zijn wijze van ontstaan om,
nieuwe gebouwen zeer eenvoudig aan te brengen,
als elementen waarvan de elementen van het oude of systemen daarvan worden genomen,
in nieuwe rangschikking,
waarbij men de rangschikking in het oude gebouw voor oogen behoudt.
Op de mogelijkheid van zulk bouwen van nieuwe systemen in bepaalden samenhang
met een vooraf gegeven systeem,
komt neer,
wat men noemt de `eigenschappen' van het gegeven systeem.}

\dutch{En een belangrijke rol speelt bij den opbouw der wiskunde
juist dat \emph{inpassen in een gegeven systeem} van nieuwe systemen,
dikwijls in den vorm van een onderzoek
naar de mogelijkheid of onmogelijkheid van een inpassing,
die aan bepaalde voorwaarden voldoet,
en in geval van mogelijkheid naar de verschillende wijzen waarop.}
\citep[pp.~77–78, original emphasis]{Brouwer1907}
\end{quotation}}
\end{quotation}

Brouwer calls the whole of the old building,
the new one,
and the correspondences between them,
a `fitting-in'
(Dutch: \dutch{inpassing}).
There is a familiar act-object ambiguity here.
A fitting-in,
as an object,
is given by two buildings and a correspondence between
elements of these buildings.
A fitting-in as an act is the act of building up a
fitting-in as an object.
In either sense,
the concept of property in Brouwer's sense stands in
contrast
to a primarily
\emph{logical}
(and hence linguistic)
concept of property.

In my view,
a fitting-in,
as an object,
is what is elsewhere known
as a state of affairs
(\german{Sachverhalt}).
This brings me to
Parsons' observation that
\begin{quote}
Brouwer is not as clear as he might be about the
distinction between intuition
\emph{of}
and intuition
\emph{that}.
Writers about Brouwer
tend to be even less so.
\citep[p.~175n71]{Parsons2008} 
\end{quote}
I confess that I,
as a writer about Brouwer,
have never stated my view on this clearly;
but it has always been that intuition
\emph{that}
is a special case of intuition
\emph{of},
namely,
of a state of affairs.
In this I follow Husserl,
who had learned the term
`\german{Sachverhalt}',
in this particular use,
from Carl Stumpf,
his teacher in Halle;
and the concept itself had been used by
the teacher of both,
Franz Brentano
\citep{Smith1989}.
The fact that Husserl saw things this way
is remarked on by Parsons:
\begin{quote}
The basic notion for Husserl is intuition of;
nonetheless intuition is what distinguishes
different ways of entertaining a proposition
from actively knowing it.
But it seems to me that Husserl reduces
intuition that to a form of intuition of,
where the object is not what we would call
a proposition but rather a state of affairs
(Sachverhalt).
\citep[143n10]{Parsons2008}
\end{quote}
and
\begin{quote}
Husserl seems to regard intuition
\emph{that}
as a species of intuition
\emph{of}:
Evidence of a judgement is a situation in which
the state of affairs that obtains if it is true is
`itself given'.
Because,
typically,
a proposition involves reference to objects,
evidence will involve intuition of those objects,
but they play the rôle of constituents of a state of affairs
that is also intuitively present,
at least in the ideal case.
\citep[146]{Parsons2008} 
\end{quote}
In a footnote,
Parsons adds:
\begin{quote}
In his discussion of truth,
Husserl talks about the `ideal of final fulfillment'
(LU VI [\citet{Husserl1984b}] §§37–39).
Later,
he concedes that this is in interesting cases
not achieved or even achievable,
so that final fulfilment is a kind of Kantian idea.
\citep[146n21]{Parsons2008} 
\end{quote}
But it can be argued that,
while this is true in general,
Husserl's conception
of a particular class of objects,
the purely categorial objects,
does not in fact allow 
the use of Kantian ideas
to conceive of the fulfilment of
intentions directed towards them.
That class includes the objects of pure mathematics;
for further discussion,
I refer to \citet[pp.~78–79]{Atten2010}.

Let me now return to induction
and state the assertion-condition for 
$\forall x(A(x) \rightarrow A(Sx))$ 
as follows:
I should have a construction method $f$ that,
given a constructed number object $x$,
yields a construction method $g=f(x)$ to construct a
fitting-in $A(Sx)$
whenever I am given a fitting-in $A(x)$,
so that $g(A(x))=A(Sx)$.%
\footnote{Note that here `$A(Sx)$' stands for a mathematical
construction in intuition, 
not for a proposition.
The ambiguity of the notation will be resolved by the context.}
But if that condition is fulfilled,
I can combine these two methods $f$ and $g$,
together with the device of ordered pairing,
into
one method $h$ that,
given an ordered pair of construction objects $x$ and the fitting-in
$A(x)$,
yields the ordered pair of construction objects
$Sx$ and the fitting-in $A(Sx)$:
Define
\begin{equation*}
h(\langle x,A(x)\rangle)=\langle Sx,g(A(x))\rangle
\end{equation*}
Then
\begin{equation*}
h(\langle x,A(x)\rangle)=\langle Sx,A(Sx)\rangle
\end{equation*}
As
$h$
is uniform in its operation on the two components of the pair,
by our earlier consideration we have intuitive knowledge
that we can iterate it.
Given a construction for $0$ and,
by hypothesis,
a construction for $A(0)$,
we construct in intuition the ordered pair
$\langle 0,A(0)\rangle$ and iterate $h$:
\[
\begin{array}{lcl}
\langle 0,A(0)\rangle & & \\
h(\langle 0,A(0)\rangle) & = & \langle 1,A(1)\rangle \\
h(\langle 1,A(1)\rangle) & = & \langle 2,A(2)\rangle\\
 & \vdots & \\
\end{array}
\]
Thus we have a uniform construction
of fitting-ins of
$n$
into
$A(n)$,
justifying the conclusion of the induction principle.
Again,
on Brouwer's account
not just the individual members of this
potentially infinite sequence are given in intuition,
but also that sequence as such.

The generality of the account,
and hence of the validity of induction,
with respect to the predicate $A$
follows because it imposes no conditions on the predicate
other than that it be total
and that it be constructive in Brouwer's
sense.
The generality claim
is not grounded on 
a prior overview over the domain of
all constructible predicates,
but on general knowledge I have,
through reflection on acts,
of the genesis of fitting-ins.
Casting a mathematical view on acts of mathematical construction
that proceed iteratively,
one identifies their common structure and thereby obtains
iteration as an act type,
itself given in intuition.
It is in this sense that the judgement that the 
induction principle is valid is the result
of second-order mathematics.

Parsons writes that
\begin{quote}
What [Brouwer] calls the
`original intuition of mathematics'
is not an intuition
\emph{of}
iteration or of the natural numbers.
I think one can regard Brouwer as holding that any natural
number can
be given in intuition;
iteration and the structure of the natural numbers arise
through the `self-unfolding' of the intuition,
but there is no reason to suppose that either is an
\emph{object}
of intuition.
The phrase
`intuition of iteration'
does not,
so far as I know,
occur in Brouwer's writings.
\citep[p.~214]{Parsons1986}%
\footnote{The quotation continues:
`it
\emph{was} used by Hermann Weyl,
who said that on the
basis of the intuition of iteration we are convinced that
the concept
of natural number is ``extensionally definite'
(\emph{umfangsdefinit}; \citet[p.~85]{Weyl1919}),
that is,
that the natural numbers are a domain over which
classical quantification is valid.
In my opinion,
such a view at best presupposes a different conception of
intuition
[than Kant's or Brouwer's,
it seems,
given the following sentence]
and is at worst confused.
In fact,
Weyl's conception of intuition seems to derive not from Kant
or Brouwer but from Husserl.'}
\end{quote}
However,
as we have seen,
a phrase that does occur in Brouwer's writings is that `and so on' is
`immediately read off from the Urintuition'
[1907:180];
and,
in a handwritten note to that passage,
that
`and so on' is among the `polarizations of the Urintuition'
\citep[p.~136n108]{Dalen2001a}.
These formulations imply that `and so on' is
given as part of the Urintuition,
understood in its self-unfolding,
as `and so on' involves more than one act;
specifically,
as what Husserl would call a dependent part,
a part that cannot be given independently.
If this is combined with Brouwer's understanding of 
`and so on'
as
`the indefinite repetition
of
one and the same
object or operation,
even if that object or operation is defined in a rather complex way',
which implies generality,
then it seems clear that Brouwer
does think of iteration as an
object of intuition.
As reconstructed here,
it is given as an act type,
the intuition of which is founded on
objectified iterative acts,
given in intuition by reflection.
Thus,
Brouwer disagrees with Parsons here,
and also with Tait,
who both hold that the general idea of iteration
is `not found in intuition'
\citep[p.~539]{Tait1981},
\citep[p.~225]{Parsons1986}.%
\footnote{Tait's rejection of intuition goes even further:
\begin{quote}
However and in whatever sense
one can represent the operation of successor,
to understand Number
one must understand the idea of iterating this operation.
But to have this
idea,
itself not found in intuition,
is to have the idea of Number
\emph{independent of any sort of representation in
intuition}.
The
same objection applied to Brouwer's analysis of number in
terms of
consciousness of succession in time (two-ity).
Again,
the essence of
number is in the iteration of that operation,
and the idea of iteration is
not founded on time consciousness.
\citep[p.~539–540]{Tait1981}
\end{quote}
I agree with Parsons' comment
on that passage:
\begin{quote}
\label{ParsonsOnTait}
[This] seems to imply either that the `idea of number'
is a concept of an abstract structure that does not depend
on any manner in which an instance of the structure might be
given,
or that an instance is given in an essentially non-intuitive
way.
Tait does not argue for either of these positions,
and I am inclined to reject both.
\citep[p.~225n16]{Parsons1986}
\end{quote}
}

\section{Some questions to Parsons}
\label{QuestionsParsons}
The account in the previous section suggests various reflections
on Parsons' arguments.

As Parsons notes,
\begin{quote}
Although the concept of a string of strokes involves
iteration,
the proposition that every such string can be extended is
not
an inductive conclusion. 
A proof by induction would be
circular.
\citep[p.~175]{Parsons2008}
\end{quote}
In the Brouwerian account,
the circularity was,
in effect,
blocked by arguing that there is one
iterative
form that is given to us without our having to construct it
first –
the structure of inner time awareness.
Parsons,
in his own setting,
proposes a different type of solution,
which is to say that,
just to
\emph{know}
that a string of strokes can be extended,
we do not have to think of the string we are extending as
having
been obtained by iterated application of adding one more
string.
We can make do, 
Parsons says \citeyearpar[p.~175]{Parsons2008},
with a `proto-conception' of string,
in which we so to speak willingly forget about those
iterations and
then add one more stroke to it.
The idea is that the proto-conception is rich enough
to make us see extendability of every string,
but too poor to set a circularity in
motion.
But here one could,
I think,
ask a question about type-preservation again.
How do we know that the result is a
\emph{string}?
Parsons explains that
\begin{quote}
to see the
\emph{possibility}
of adding one more,
it is only the general structure that we use,
and not the specific fact that what we have before us
was obtained by iterated additions of one more.
This is shown by the fact that,
in the same sense in which a new stroke can be added
to any string of strokes,
a new stroke can be added to any bounded geometric
configuration.
\citep[p.~175]{Parsons2008}
\end{quote}
But precisely because of that generality
– the independence of the type of the given geometric
configuration – 
it is not clear to me that this argument
gives us a purchase
on the
type
of the particular resulting configuration.%
\footnote{A good discussion and reply to earlier criticisms
of the idea that,
on the conception of intuition set out by Parsons,
it is intuitive knowledge that any string can be extended by
one stroke,
is \citet[section 4]{Jeshion2014}. (Thanks to Charles
Parsons for drawing my attention to this paper.)}

Parsons acknowledges 
that the
limits of intuitive
knowledge that he arrives at are rather narrow,
and that many will hold that this is due to the very
restricted character
of the conception of intuition he develops 
\citep[p.~316]{Parsons2008}.
One particularly strong constraint Parsons works with is
Kantian in the broad sense that intuition is
intuition of spatio-temporal objects.
However,
Parsons explicitly leaves open that there might be different
models of
intuition,
on which there would be intuition also of other types of
objects,
and he mentions Husserl and Gödel.
I think that material in the Gödel archive,
reading notes and work for the revision of the Dialectica
Interpretation,
shows that he was much more committed to Husserl's notion
than he was willing to let on in either his publications
or his conversations with people
without too much interest in phenomenology 
\citep[section 1.2 and chs 4, 6]{Atten2015a}.
Be that as it may,
the following two questions to Parsons are
directly concerned with Husserl's notion.

As we saw in section \ref{Parsons},
Parsons says that
we have intuitive knowledge
\emph{that}
$p$ if $p$ can be `seen' to be true
on the basis of intuiting the objects that $p$ is about.
But he does not say much more about how we see this.
Does the `seeing' involve any intuitive component other
than that of the objects?
We would arrive at something like Husserl's 
`categorial intuition',
but I am not sure Parsons would be willing to embrace that.%
\footnote{In
`Mathematical intuition',
Parsons remarks that
`Husserl does undertake to show that in categorial intuition
there is something analogous to sensations in sense perception.
In my view,
he lapses into obscurity in explaining this
(LU, VI [\citet{Husserl1984b}], §56).
I am not sure to what extent this can be cleared up'
\citep[p.~166n8]{Parsons1980}.
In that 
section,
Husserl outlines a theory
of the `categorial representation'
(`\german{kategoriale Räpresentant}'),
which,
however,
he later dropped \citep[p.~535]{Husserl1984b}, \citep{Lohmar1990}.
Note that Parsons' critical remark about it
has not been included in 
\emph{Mathematical Thought and Its Objects}.
An alternative interpretation of categorial
intuition
is proposed in \citet{Atten2015b}.}
But if,
alternatively,
the component that the seeing involves beyond intuition
of the objects is itself non-intuitive,
what would justify calling the resulting knowledge
`intuitive knowledge'?

The following question is motivated by the fact
that,
although Parsons mentions Husserl in
\emph{Mathematical Thought and Its Objects},
he does not mention the aspect of Husserl
that brings Husserl closest to Brouwer,
inner time awareness.
It would seem to be a natural question,
however,
whether Parsons would not be willing to extend his notion
of intuition,
such as he explicitly describes it,
with a Brouwerian or Husserlian intuition of inner time.
It is not clear to me whether Parsons has a principled
reason not to.
One may of course observe that,
given his long-standing engagement with both Brouwer's and
Husserl's
thought,
if he had wanted to exploit their notion of intuition of
inner time,
he would have done so by now.

Yet in
\emph{Mathematical Thought and Its Objects},
Parsons appeals to Brouwer's intuition of two-ity twice.

First,
as one way of seeing that every string of strokes can be
extended:
\begin{quote}
Let us return to the proposition that any string can be
extended.
The idea that this rests on a capability of the mind is a
very natural one and in certain respects acceptable.
I have proposed two different ways of seeing this,
one resting on the figure-ground structure of perception and
one (Brouwer's) resting on temporal experience.
… We experience the world as temporal,
and have the conviction that we can continue into a
proximate future,
in which the immediate past is retained 
\citep[p.~177]{Parsons2008}.%
\footnote{Also in `Intuition in constructive mathematics':
\begin{quote}
The `original intellectual phenomenon of the falling apart
of a life moment
into two qualitatively distinct things' 
\citep[p.~153]{Brouwer1929A},
is in a certain way iterable:
since we can divide our experience into past and
present/future,
independently of its objects,
we can continue to repeat that division, so
that there `arises by self-unfolding of the original
intellectual phenomenon
the temporal series of appearances of abitrary multiplicity'
(ibid.)'
\citep[pp.~212\textendash{}213]{Parsons1986}
\end{quote}}
\end{quote}

The second appeal is made in a
in a comment on Bernays.
Bernays had written:
\begin{quote}
We are conscious of the freedom we have to advance from one
position arrived at in the process of counting to the next
one.%
\footnote{[\german{Zunächst sind wir uns der Freiheit bewusst,
von einer erreichten Stelle
im Zählprozess jweils noch um Eins
fortzuschreiten} \citep[p.~469]{Bernays1955}]}
\end{quote}
Parsons comments:
\begin{quote}
It would take some argument to show that there is no appeal
here to the temporal character of experience,
such as we find in Brouwer.
\citep[p.~337]{Parsons2008}
\end{quote}
As we just saw,
Parsons is willing to appeal to Brouwer's temporal
experience
as one way of upholding the claim that every string of
strokes can be
extended;
and the present comment on Bernays
is made by Parsons in the context of defending his own
claim that purely rational evidence cannot replace appeals
to intuition completely.
I therefore take also this comment to show that
Parsons indeed wishes to accept
Brouwer's temporal intuition.

It is,
of course,
characteristic of Parsons' general approach 
in the philosophy of mathematics that,
to convince us that
we have the capability to extend any string,
he refers to two types of
\emph{intuition}
– perception
and temporal experience.
However,
in his description and use of Brouwer's two-ity,
it is the discrete elements in the two-ity and their
order that he exploits,
never quite mentioning the continuum in between the two
things.
Intuitionists subscribe to the view that
the intuitions of the discrete and the continuous
cannot be accepted independently from one another 
\citep[p.~8]{Brouwer1907}/\citep[p.~17]{Brouwer1975}.
A question to Parsons,
then,
is whether he means to do just that,
and if so,
whether that is possible.

I should like to close with a remark on impredicativity.

Parsons 
\citeyearpar[pp.~293\textendash{}294]{Parsons2008}
accepts the following argument,
proposed by
Dummett,
to the effect
that the notion of natural number is impredicative:
\begin{quote}
The totality of natural numbers is characterised as one
for which induction is valid with respect to any
well-defined property,
where by a `well-defined property' is understood one which
is well-defined
relative to the totality of natural numbers.
In the formal system,
this characterisation is of course weakened to
`any property definable within the formal language';
but the impredicativity remains,
since the definitions of the properties
may contain quantifiers whose variables 
range over the totality characterised.
\citep[p.~199]{Dummett1978}
\end{quote}

It is curious that Dummett,
who raised the issue in his paper
on Gödel's incompleteness theorem of 1963,
does not discuss it in his later book
\emph{Elements of intuitionism} \citep{Dummett1977,Dummett2000b}.
Be that as it may,
John Myhill has argued that the constructivist
can avoid this problem by saying that
to the constructivist,
the notion `finite' or some equivalent idea
such as `natural number' or `ancestral' is clear whereas
impredicative
definitions are not
\citep[p.~27]{Myhill1974}.
Parsons 
answers to Myhill
that that reply depends
on `a dogmatic view of the clarity of the notion of
natural number and the evidence of mathematical induction'.
He adds that `such a dogmatic view could plausibly be
attributed to Poincaré and possibly also Brouwer' \citep[p.~294n44]{Parsons2008}.

I have tried to show that
the view is,
in Brouwer's setting,
less dogmatic than it may seem.
Moreover,
as we saw above
(p.~\pageref{BrouwerProperty}),
for Brouwer a (mathematical) property is primordially
not a logically defined predicate,
but
a fitting of one
mathematical building into another.
As the notion of mathematical building
depends on that of the two-ity and its unfolding,
this means that
the notion of
property presupposes an intuition that by
itself
suffices to
give the natural numbers.%
\footnote{Even if one conceives of properties in
terms of logically defined predicates,
the intuition of `and so on' is presupposed,
namely in our grasp of the syntax.}
This contradicts Parsons' conclusion that
`the concept of natural number cannot determine what
counts as a well-defined predicate'
\citep[p.~267]{Parsons2008};
but Brouwer and Parsons are speaking from different backgrounds.

\paragraph{Acknowledgement.} 
I express my profound gratitude to Charles Parsons,
for his generous instruction and friendship
over the years.

This is a revised and expanded version of a paper presented 
in Jerusalem 
on December
4,
2013,
at
the conference on the work of Parsons,
`Intuition and Reason',
held in Tel Aviv and Jerusalem,
Dec 2–5,
2013.
Earlier versions were
presented at `\french{Vuillemin,
lecteur de Kant}' at the Archives
Poincaré in Nancy
(Dec 15,
2012);
at the seminar
`\french{Mathématiques et Philosophie,
19e et 20e siècles}' at SPHERE,
Paris
(Feb 21,
2013);
at the seminar `\french{Les usages de la
phénoménologie
(I):
temps vécu et temps cosmique}' in Lille
(Mar 28,
2013);
and at the joint PHILMATH seminar of SND and IHPST,
Paris
(Oct
28,
2013).
I thank the organisers for their invitations,
and the
audiences and later readers for their questions and comments;
in particular,
in addition to Parsons:
Dirk van Dalen,
Gerhard Heinzmann,
Wilfried Sieg,
Göran
Sundholm,
Bill Tait,
Joseph Vidal-Rosset,
Albert Visser,
and Michael Wright.

\paragraph{Postscript, May 2015}

When this paper had almost been completed,
Charles Parsons
shared with me his very recent manuscript `Intuition
revisited'.
Among other things,
it contains a revision of
the conception of intuition developed in chapter 5 of
\emph{Mathematical Thought and Its Objects},
in the light of
criticism by Felix Mühlhölzer
(Mühlhölzer 2010).
Parsons
revises that conception by replacing spatial intuition with
temporal intuition as the foundation of intuitive knowledge
in arithmetic,
with reference to Brouwer.
The question
whether this revision has a bearing on the evidence of
induction and recursion is explicitly left aside.
A
remaining difference between Parsons and Brouwer is that in
Parsons' exposition,
external objects still have a role,
so
as to avoid charges of subjectivism or solipsism.
A
Brouwerian reply,
based on a Husserlian reading of Brouwer's
writings and correspondence,
would be to say that the
structure of inner time awareness is identical across
different minds because it is a structure of
(not empirical
but)
transcendental subjectivity
(Van Atten 2004,
chapter
6).
But Parsons does not want to commit himself to a
transcendental view of consciousness.

\bibliography{Induction}

\end{document}